\documentclass[reqno,a4paper,12pt]{article}
\usepackage{hyperref}
\usepackage{amsmath,amssymb,amsthm,authblk}
\usepackage{textcomp}
\usepackage[pdftex]{graphicx}
\usepackage{amsmath,epsfig}
\usepackage{breqn}
\usepackage{epstopdf}
\usepackage{epsfig}
\usepackage{latexsym}
\usepackage{xcolor}
\usepackage{mhsetup}
\usepackage{mathtools}
\usepackage{mathrsfs}
\usepackage{float}
\usepackage{xfrac}
\usepackage[english]{babel}

\def\r{\rho}

\def\vp{\varphi}

\usepackage{fancyhdr}

\usepackage{titlesec}

\title{\bf\large draft 0}
\date{\vspace{-5ex}}

\author {\small Shadi Al-Omari}
\affil{\small\it The Preparatory Year Program, King Fahd
University of Petroleum and Minerals, Dhahran 31261, Saudi
Arabia.\\
\small\it The Interdisciplinary Research
Center in Construction and Building Materials, Dhahran 31261, Saudi
Arabia. } \affil{shalomari@kfupm.edu.sa}

\begin{document}
\title{\bf\large New Decay Results for a Partially Dissipative Viscoelastic Timoshenko System with Infinite Memory}\maketitle

\begin{abstract}
In this paper, we consider the following dissipative viscoelastic with memory-type Timoshenko system
\begin{equation*}
\begin{gathered}
\begin{cases}
\rho_1 \phi_{tt}  - \kappa ( \phi _{x}  + \psi) _x + \kappa \int_0^\infty g(s) (\phi_x +\psi)_x(t-s) ~ds =0   & \text{in}~ \left( {0,L} \right) \times \mathbb{R}^+ , \\
\rho_2 \psi_{tt} - b \psi_{xx} + \kappa ( \phi _{x}  + \psi)-\kappa \int_0^\infty g(s) (\phi_x +\psi)(t-s)~ ds=0 & \text{in}~ \left( {0,L} \right) \times \mathbb{R}^+ , \\
\end{cases}
\end{gathered}
\end{equation*}
with Dirichlet boundary conditions, where $g$ is a positive
non-increasing function satisfying, for some nonnegative functions
$\xi$ and $H$, \[g'(t)\leq-\xi(t)H(g(t)),\qquad\forall~ t\geq0.\]
Under appropriate conditions on $\xi$ and $H$, we establish some new decay results for the case of equal-speeds of propagation that generalize and improve many earlier results in the literature.
\end{abstract}
\hspace{-0.17 in}\bf{Keywords:}\rm \;  Timoshenko system, Infinite memory, General decay, Convex functions, equal wave speeds .\\
\bf{AMS Classification:}\rm \; $35\text{B}37\cdot 35\text{L}55
\cdot 74\text{D}05 \cdot 93\text{D}15 \cdot 93\text{D}2$.

\numberwithin{equation}{section}
\newtheorem{theorem}{Theorem}[section]
\newtheorem{lemma}[theorem]{Lemma}
\newtheorem{remark}[theorem]{Remark}
\newtheorem{corollary}{Corollary}[section]
\newtheorem{definition}[theorem]{Definition}
\newtheorem{proposition}[theorem]{Proposition}
\newtheorem{example}{Example}
\newtheorem{Corollary}[theorem]{Corollary}
\allowdisplaybreaks
\section{Introduction.}
We consider the following Timoshenko system with a viscoelastic
dissipation mechanism coupled on the shear force:
\begin{equation}\label{1.1}
\begin{cases}
\begin{array}{ll}
\rho_1 \phi_{tt}  - \kappa ( \phi _{x}  + \psi) _x + \kappa \int_0^\infty g(s) (\phi_x +\psi)_x(t-s) ~ds =0   & \text{in}~ \left( {0,L} \right) \times \mathbb{R}^+ , \\
\rho_2 \psi_{tt} - b \psi_{xx} + \kappa ( \phi _{x}  + \psi)-\kappa \int_0^\infty g(s) (\phi_x +\psi)(t-s)~ ds=0 & \text{in}~ \left( {0,L} \right) \times \mathbb{R}^+ , \\
\phi (0, t) = \phi (L, t) = \psi_x (0, t) = \psi_x (L, t) = 0,& t \geq 0,\\
\phi_x (x, -t) = \phi_0(x,t),~ \phi_t(x, 0) = \phi_1(x),& x \in  (0,L),\\
\psi (x, -t) = \psi_0(x,t),~ \psi_t(x, 0) = \psi_1(x),& x \in  (0,L),\\
\end{array}
\end{cases} \tag{$P$}
\end{equation}
where $(x, t) \in (0, L) \times (0, \infty)$, $\rho_1$, $\rho_2$, $b$ and $\kappa$ are positive constants, $L > 0$ is the length of the beam and $\mathbb{R}^+ = (0,\infty )$. Corresponding to the unknown variables $\phi$  and $\psi$,  and $\phi_0$, $\phi_1$, $\psi_0$ and $\psi_1$ as an initial data and $g$ is the relaxation function which is also known as memory kernel satisfying some conditions to be specified in the next section. The Timoshenko model was first derived in 1921 by Timoshenko \cite{Timoshenko1921, timoshenko1922} to describe the dynamics of a beam by taking the transverse shear strain into consideration.\\
Giorgi et al. \cite{Giorgi} considered the following semilinear
hyperbolic equation with linear memory in a bounded domain $\Omega
\subset \mathbb{R}^3$
\begin{equation}\label{giorgi}
    u_{tt}- K(0)\Delta u -\int_{0}^{+\infty}K^{\prime}(s) \Delta u(t-s)ds+g(u)=f\hspace{0.1in}\text{in}\hspace{0.05in}\Omega\times
    \mathbb{R}_+,
\end{equation}
with $K(0), K(+\infty)>0$ and $K^{\prime}\le 0$ and proved the
existence of global attractors for the solutions. Conti and Pata
\cite{Pata2} considered the following semilinear hyperbolic
equation:
\begin{equation}\label{E:Pata}
u_{tt}+\alpha u_{t}- K(0)\Delta u -\int_{0}^{+\infty}K^{\prime}(s)
\Delta
u(t-s)ds+g(u)=f\hspace{0.1in}\text{in}\hspace{0.05in}\Omega\times\mathbb{R}_+,
\end{equation}
where the memory kernel is a convex decreasing
smooth function such that $K(0)>K(+\infty)>0$ and
$g:\mathbb{R_+}\rightarrow \mathbb{R_+}$ is a nonlinear term of at
most cubic growth satisfying some conditions. They proved the
existence of a regular global attractor. In \cite{Appleby},
Appleby et al. studied the linear integro-differential equation
\begin{equation}\label{appleby}
    u_{tt}+Au(t)+\int_{-\infty}^{t}K(t-s)Au(s)ds=0 \hspace{0.1in}\text{for} \hspace{0.1in}
    t>0,
\end{equation}
and established an exponential decay result for strong solutions
in a Hilbert space. Pata \cite{pata4} discussed the decay
properties of the semigroup generated by the following equation:
\begin{equation}\label{pata}
    u_{tt}+\alpha Au(t)+\beta u_{t}(t)-\int_{0}^{+\infty}\mu(s)Au(t-s)ds=0  \hspace{0.1in}\text{for} \hspace{0.1in}
    t>0,
\end{equation}
where $A$ is a strictly positive self-adjoint linear operator and
$\alpha>0, \beta\ge 0$ and the memory kernel $\mu$ is a decreasing
function satisfying specific conditions. Subsequently, they
established necessary as well as the sufficient conditions for the
exponential stability. In \cite{Guesmia22}, Guesmia considered
\begin{equation}
u_{tt}+Au-\int_{0}^{+\infty}g(s)Bu(t-s)ds=0
\hspace{0.1in}\text{for} \hspace{0.1in}
    t>0,
\end{equation}
and introduced a new ingenuous approach for proving a more general
decay result based on the properties of convex functions and the
use of the generalized Young inequality. He used  a larger class
of infinite history kernels satisfies the following condition
\begin{equation}\label{gg}
\int_{0}^{+\infty}\frac{g(s)}{H^{-1}(-g^{\prime}(s))}ds+\sup_{s\in\mathbb{R}_+}\frac{g(s)}{H^{-1}(-g^{\prime}(s))}<+\infty,
\end{equation}
such that
\begin{equation}\label{E:3}
H(0)=H^{\prime}(0)=0\hspace{0.05in}\text{and}\hspace{0.05in}\lim_{t\rightarrow
+\infty}{H^{\prime}(t)}=+\infty,
\end{equation}
where $H:\mathbb{R}_+\to \mathbb{R}_+$  is an increasing strictly
convex function. Using this approach, Guesmia and Messaoudi
\cite{guesmia2} later looked into
\begin{equation*}
u_{tt}-\Delta u + \int_{0}^{t}g_{1}(t-s) div (a_{1}(x) \nabla
u(s))ds + \int_{0}^{+\infty}g_{2}(s)div (a_{2}(x) \nabla
u(t-s))ds=0,
\end{equation*}
in a bounded domain and under suitable conditions on $a_{1} \;
\text{and} \; a_{2}$ and for a wide class of relaxation functions
$g_{1}\; \text{and} \; g_{2}\;$ that are not necessarily decaying
polynomially or exponentially and established a general decay result
from which the usual exponential and polynomial decay rates are only
special cases. Recently, Al-Mahdi \cite{al2020stability} consider
the following viscoelastic plate problem with a velocity-dependent
material density and a logarithmic nonlinearity:
\begin{equation}\label{First}
\vert u_t \vert^{\rho} u_{tt}+\Delta^2 u+\Delta^2
u_{tt}-\int_{0}^{+\infty}g(s)\Delta^2 u(t-s)ds=k u \ln{\vert u\vert}
\hspace{0.2in}\text{in   } \Omega\times (0,\infty),
\end{equation} where $\Omega $ is a bounded domain of $\mathbb{R}^2,$ with a
smooth boundary $\partial \Omega$. He established an explicit and
general decay rate results with imposing a minimal condition on the
relaxation function, that is,
\begin{equation}\label{kernel}
g^{\prime}(t)\le -\xi(t) H(g(t)),
\end{equation}where the two functions $\xi$ and $H$  satisfy some
conditions. Very recently, Al-Mahdi \cite{al2020general} considered
the following plate problem:
\begin{equation*}
u_{tt}-\sigma \Delta u_{tt}+ \Delta^2
u-\int_{0}^{+\infty}g(s)\Delta^2 u(t-s)ds=0,
\end{equation*}
and proved that the stability of this problem holds for which the
relaxation function $g$ satisfies the condition \eqref{kernel}\\\\
For Timoshenko systems with infinite memory, Rivera et al.
\cite{Rivera} considered vibrating systems of Timoshenko type with
past history acting only in one equation. They showed that the
dissipation given by the history term is strong enough to produce
exponential stability if and only if the equations have the same
wave speeds. In the case that the wave speeds of the equations are
different, they showed that the solution decays polynomially to zero
if the corresponding system does not decay exponentially as time
goes to infinity, with rates that can be improved depending on the
regularity of the initial data. Guesmia et al. \cite{Guesmia2012}
have adopted the method introduced in \cite{Guesmia22} with some
necessary modifications to establish a general decay of the solution
for a vibrating system of Timoshenko type in a one-dimensional
bounded domain with an infinite history acting in the equation of
the rotation angle. Guesmia and Messaoudi \cite{Guesmiaetal3}
discussed a Timoshenko system in the presence of an infinite memory,
where the relaxation function satisfies
$g'(t)\leq-\xi(t)g(t),~\forall t\geq0$ and established some general
decay results for the equal and nonequal speed propagation cases.
Recently, Guesmia \cite{Guesmiarecent} adapted the approach of
\cite{Mus-2017} to two models of wave equations with infinite memory
and proved, under the condition  $g'(t)\leq-\xi(t)G(g(t)),~\forall
t\geq0$, where $\xi$ is satisfying $\int_0^{+\infty} \xi(s) ds = +
\infty$ decay rate of solutions and the growth of $g$ at infinity.
Al-Mahdi \cite{al2020stability,al2020general} also adapted the
approach of \cite{Mus-2017} to some viscoelastic plate equations
with relaxation functions satisfy the condition
$g'(t)\leq-\xi(t)G(g(t)),~\forall t\geq0$. The results of
\cite{Guesmiarecent} and \cite{al2020stability,al2020general}
improved and generalized the ones of \cite{Guesmia1},
\cite{Guesmia22},\cite{Guesmia-abstract}, \cite{youkana_lemma1_a},
\cite{almahdi2} and \cite{almahdibvp} by getting a better decay rate
and deleted some assumptions on the boundedness of initial data.

The rest of this paper is organized as follows. In section 2, we
present some assumptions and material needed for our work. Some
technical lemmas are presented and proved in section 3. Finally,
we state and prove our main decay results and provide some
examples in section 4.\\

\section{Preliminaries} In this section, we introduce some notation and assumption, present some useful lemmas and state the existence theorem. Let us start by introducing the following standard functional spaces:
\begin{equation*}
\begin{array}{lc}
L^2 := L^2(0,L),~~~ ||u||_2^2 =\int_0^L | u(x)|^2 dx, \\
H^1 := H^1(0,L),~~~ ||u||^2_{H^1} = ||u_x||_2^2+ ||u||_2^2, \\
L^2_{\ast} := L^2_{\ast} (0,L) =\Big\{ u \in  L^2(0,L) ;\frac{1}{L}\int^L_0 u(x) dx = 0 \Big\} , \\
H^1_0 := H^1_0 (0,L) =\Big\{ u \in  H^1(0,L) ; u(0) = u(L) = 0 \Big\}, \\
H^1_{\ast}  := H^1_{\ast}  (0,L) = \Big\{ u \in  H^1(0,L) ;\frac{1}{L}\int^L_0 u(x) dx = 0 \Big\}. \\
\end{array}
\end{equation*}
Due to Poincar\'e's inequality, we can also consider the equivalent norms in $H^1_0$ and $H^1_\ast$,
\begin{equation*}
 ||u||_{H^1_0}= || u_x||_2  ~~~~~~\text{and}~~~~~~ || u||_ {H^1_\ast}= ||u_x||_2,
\end{equation*}
respectively. In this work, we will always denote by $c_p > 0$ the Poincar\'e constant.\\
\textit{\textbf{Assumptions:}} We assume that the relaxation function $g$ satisfy the following hypotheses.\\
\textit{\textbf{(A1)}} $g : \mathbb{R}^+ \rightarrow  \mathbb{R}^+$ is a non-increasing differentiable function such that
\begin{equation} \label{2.5}
g(0) > 0 ~~~~~~\text{and}~~~~~~ \ell := 1-\int^\infty_0 g(s)ds > 0.
\end{equation}
and
\begin{equation} \label{2.51}
\int^\infty_0 g(s)ds > max \biggl\{ 
\frac{31}{32}, \frac{64 \rho_1 L^2}{64 \rho_1 L^2 + \rho_2}
\biggr\}.
\end{equation}
\textit{\textbf{(A2)}} There exist a non-increasing differentiable
function $\xi:\mathbb{R}^+ \rightarrow  \mathbb{R}^+$ and a
$C^1-$function $H:[0, +\infty) \rightarrow [0, +\infty)$ which is
either linear or strictly increasing and strictly convex
$C^2$-function on $(0,r]$ for some $r > 0$ with
 $H(0)=H^{\prime}(0)=0$,~$\lim_{s\rightarrow + \infty}
H^{\prime}(s)=+\infty $,~$s\mapsto s H^{\prime}(s) $ and $s \mapsto
s \left(H^{\prime}\right)^{-1}(s)$ are convex on $(0,r]$. Moreover,
there exists a positive nonincreasing differentiable function $\xi$
such that
\begin{equation}\label{2.5.1}
g^\prime(t) \leq - \xi (t) H(g(t)) , ~~~~ \forall  ~~ t \geq 0.
\end{equation}

\begin{remark}\label{new1}
The condition (\ref{2.51}) means that the area under the graph of  $g$ is bounded below. It does not contradict the relation (\ref{2.5}) since
$$ C_0 := max \biggl\{ \frac{31}{32},\frac{64 \rho_1 L^2}{
64 \rho_1 L^2 + \rho_2} \biggr\} < 1.$$
\end{remark}
\begin{remark}\label{new2}
From (\ref{2.51}) we infer that there exists
a time $t_0 > 0$ large enough such that 
\begin{equation}
1- h(t) = \int^t_0 g(s) ds \geq \int^{t_0}_0 g(s) ds := g_0 > C_0~~~~ \forall  t \geq  t_0.
\end{equation}
\end{remark}


\begin{remark}\label{RG}
If $H$ is a strictly increasing, strictly convex $C^2$ function
over $(0, r]$ and satisfying  $H(0) = H'(0) = 0$, then it
has an extension $\overline{H}$, that is also strictly
increasing and strictly convex $C^2$ over $(0,\infty)$. For
example, if $H(r) = a, H'(r) = b, H''(r) = c$, and for $t
> \varepsilon$,  $\overline{H}$ can be defined by
\begin{equation}\label{hbar}
    \overline{H}(t)=\frac{c}{2} t^2+ (b-cr)t+ \left(a+\frac{c}{2} {r}^2 - b r \right).
\end{equation}For simplicity, in the rest of this paper, we use $H$ instead of $\overline{H}$
\end{remark}

\begin{remark}\label{RG2}
Since $H$ is strictly convex on $(0,r]$ and $H(0)=0,$ then
\begin{equation}\label{E2a:L2:St1}
H(\theta z)\le \theta H(z),\text{              }0\le \theta \le 1,
\hspace{0.1in}z\in (0,r].
\end{equation}
\end{remark}

\begin{theorem} Under the Assumptions (A1)-(A2) and taking $(\phi_0, \phi_1, \psi_0, \psi_1) \in  H^1_0 \times L^2 \times H^1_{\ast} \times L^2_{\ast}$, there exists a unique weak solution $(\phi , \psi )$ of problem (\ref{1.1}) in the class
$$ (\phi , \psi ) \in  C(\mathbb{R}^+;H^1_0 \times  H^1_{\ast}  ) \cap  C^1(\mathbb{R}^+;L^2 \times  L^2_{\ast} ).$$
Furthermore, if $(\phi_0, \phi_1, \psi_0, \psi_1) \in  (H^2 \cap H^1_0 )\times H^1_0 \times (H^2 \cap H^1_{\ast} )\times H^1_{\ast}$, then there exists a unique strong solution $(\phi , \psi )$ of problem (\ref{1.1}) in the class $$(\phi , \psi ) \in  C \bigl(\mathbb{R}^+; (H^2 \cap  H^1_0 ) \times  (H^2 \cap  H^1_{\ast} )\bigr) \cap  C^1(\mathbb{R}^+;H^1_0 \times  H^1_{\ast} ).$$
\end{theorem}
\begin{proof}
The proof of this theorem can be achieved by using the pattern of the Faedo-Galerkin method (see Lions book \cite{lions1969quelques}) as applied to wave equations with memory.
\end{proof}

\begin{lemma}
If $(\phi , \psi ) \in  L^2(0, T;H^1_0 \times  H^1_{\ast}), T > 0$,
then
\begin{equation} \label{2.11}
p(\cdot , t) := \phi (\cdot , t) + \widehat  \psi (\cdot , t) \in
H^1_0 (0,L).
\end{equation}
\end{lemma}
\begin{proof}
The proof follows from direct computations.
\end{proof}

We are going to see that problem (\ref{1.1}) is dissipative with only one damping mechanism given by the convolution term involving the shear force component. Indeed, under the above notation and given a weak solution $(\phi , \psi )$ of problem (\ref{1.1}), we define the corresponding energy functional\\ $ E (t) =  E \big(\phi (t), \psi (t), \phi_t(t), \psi_t(t) \big),~ t \geq  0$, by
\begin{equation} \label{2.12}
 E (t) := \frac{\rho_1}{2} ||\phi_t (t)||^2_2 + \frac{\rho_2}{2} || \psi_t (t)||^2_2 + \frac{b}{2} || \psi_x (t)||^2_2 + \frac{\kappa}{2} \gamma || p_x(t)||^2_2 + \frac{\kappa}{2} ( g \circ p_x)(t),
\end{equation}
where $\gamma:=1-\int_0^\infty g(s) ds$ and $p(x , t):= \int^x_0
u(y, t) dy$ and
\begin{equation}\label{2.8}
(g \circ u)(t) := \int_0^\infty g(s) || u(t)- u(t-s)||^2_2 ds,
\end{equation}
\begin{lemma}
The energy $ E (t)$ satisfies the following identity:
\begin{equation} \label{2.13}
\frac{d}{dt}  E (t) = - \frac{\kappa}{2} ( g^\prime \circ
p_x)(t)\leq 0,~~~~~ t > 0.
\end{equation}
\end{lemma}
\begin{proof}
Taking the multipliers $\phi_t$ and $\psi_t$ in the first two equation of \ref{1.1},
respectively, a straightforward computation leads to (\ref{2.13}).
\end{proof}

As in \cite{jin2014coupled}, we set, for any $0 < \alpha < 1$,
\begin{equation} \label{2.14}
C_{\alpha} := \int_0^{\infty} \frac{g^2(s)}{\alpha g(s) -
g^\prime(s)}ds ~~~~~~~\text{and}~~~~~~~\mu(t):= \alpha g(t)-
g^\prime (t).
\end{equation}
\begin{remark}\label{r1:4a}
Using the fact that $\frac{\alpha g^2(s)}{\alpha
g(s)-g^{\prime}(s)} <g(s)$ and recalling the Lebesgue dominated
convergence theorem, we can easily deduce that
\begin{equation}\label{E1:r1:4a}
\alpha C_{\alpha}=\int_{0}^{\infty}\frac{\alpha g^2(s)}{\alpha
g(s)-g^{\prime}(s)}ds \to 0 \text{       as      }\alpha \to 0.
\end{equation}
\end{remark}
\begin{lemma} (\cite{jin2014coupled}). Assume that assumption (A1) holds. Then for any $v \in \textbf{L}^2_{loc} (\mathbb{R}_{+};\textbf{L}^2(0,L))$, we have
\begin{equation} \label{2.15}
\int_0^L \bigg (\int_0^\infty g(s) (v(t)-v(t-s)) ds \bigg)^2 dx \leq
C_{\alpha} (\mu \circ v )(t),~~~~~ \forall ~ t \geq 0.
\end{equation}
\end{lemma}

\section{Technical lemmas} In this section, we state and proof some lemmas needed to establish our main results.
\begin{lemma}\label{cor2}
There exist a positive constant $M_0$  such that
\begin{equation}\label{h}
\begin{aligned}
\int_{t}^{+\infty}g(s)\vert \vert p_x(t)-p_x(t-s) \vert \vert_2^{2}
ds  \leq M_0 h_0(t),
\end{aligned}
\end{equation}
where $h_0(t)=\int_{0}^{+\infty}g(t+s)\left(1+\vert \vert
{p_x}_{0x}(s)\vert\vert^{2} \right) ds$
\end{lemma}
\begin{proof}
The proof of \eqref{h}  is identical to the one in
\cite{Guesmiarecent} and  \cite{al2020general}. Indeed, we have
\begin{equation}\label{hproof}
\begin{aligned}
&\int_{t}^{+\infty}g(s) \vert \vert p_x(t)-p_x(t-s)\vert \vert_2^{2}
ds \leq 2 \vert \vert p_x(t)  \vert \vert^2 \int_{t}^{+\infty}g(s)
ds\\& +2 \int_{t}^{+\infty}g(s) \vert \vert
p_x(t-s)  \vert \vert^2 ds\\
&\leq 2 \sup_{s\geq 0}\vert \vert  p_x(s)  \vert \vert^2
\int_{0}^{+\infty}g(t+s) ds +2 \int_{0}^{+\infty}g(t+s) \vert \vert
p_x(-s)  \vert \vert^2 ds\\
&\leq \frac{4 E(s)}{\kappa \gamma} \int_{0}^{+\infty}g(t+s) ds +2
\int_{0}^{+\infty}g(t+s)
\vert \vert {p_x}_{0x}(s)  \vert \vert^2 ds\\
&\leq \frac{4E(0)}{\kappa \gamma} \int_{0}^{+\infty}g(t+s) ds +2
\int_{0}^{+\infty}g(t+s) \vert \vert {p_x}_{0x}(s)  \vert \vert^2 ds\\
&\leq M_0 \int_{0}^{+\infty}g(t+s)\left(1+\vert\vert
{p_x}_{0x}(s)\vert\vert^{2} \right) ds.
\end{aligned}
\end{equation}
where $M_0=\max\big\{2,  \frac{4E(0)}{\kappa \gamma}   \big\}.$
\end{proof}

\begin{lemma}
Under assumption (A1) and (A2), the functionals $\chi_1(t), \chi_2(t)$ and $\chi_3(t)$ defined by
\begin{equation} \label{3.1}
\chi_1(t) := \int_0^L  \left( \rho_1 \phi \phi_t + \rho_2 \psi \psi_t \right) dx
\end{equation}
\begin{equation}\label{3.3}
\chi_2(t) := - \rho_1 \int_0^L \phi_t(t)  \int^\infty_0 g(s)(p(t)-p(t-s)) ds dx ,
\end{equation}
\begin{equation}  \label{3.10}
\chi_3(t) := - \rho_2 \int_0^L \psi_ t \left( h(t) p_x(t) + \int_0^\infty g(s)( p_x(t)-p_x(t-s)) ds \right) dx - \rho_2  \int_0^L \phi_t \psi_x  dx .
\end{equation}
satisfies, along with the solution of (\ref{1.1}), the estimates
\begin{equation} \label{3.2}
\chi^\prime_1(t) \leq  \rho_1 || \phi_t(t)||^2_2 + \rho_2 ||
\psi_t(t)||^2_2 - b|| \psi_x(t)||^2_2 - \frac{\kappa}{2} h(t) ||
p_x(t)||^2_2 + \frac{\kappa}{2l} C_{\alpha} (\mu \circ p_x)(t).
\end{equation}
\begin{equation} \label{3.4}
\begin{array}{ll}
\chi^\prime_2 (t) &  \leq - \frac{\rho_1 g_0}{4} || \phi_t(t)||^2_2 + \kappa (1-g_0)\frac{\rho_1 L^2}{\rho_2} h(t) || p_x(t)||^2_2 + \frac{\rho_1 L^2}{2} || \psi_t(t)||^2_2 \\
&+\frac{\rho_1 c^2_p}{g_0}   \int_0^L \left( \int_0^\infty g^\prime(s) \left(  p_x(t) - p_x (t-s) \right) ds \right)^2 dx
+ \kappa \bigg( 1 + \frac{\rho_2}{ 4(1-g_0) \rho_1 L^2}
\bigg) C_{\alpha} (\mu \circ p_x)(t).\\
\end{array}
\end{equation}
for all $t \geq  t_0$, where $c_p > 0$ is the Poincar\'e constant.
\begin{equation} \label{3.11}
\begin{aligned}
\chi^\prime_3(t) & \leq - \frac{\rho_2}{2} || \psi_t(t)||^2_2
+ 2 \kappa  [h(t)]^2 || p_x(t)||^2_2 + 2 \kappa  C_{\alpha} (\mu \circ p_x)(t) \\
& ~~~+\rho_2 g(0) g(t) || p_x(t)||^2_2
+ \rho_2 \int_0^L \left( \int_0^\infty g^{\prime}(s) (  p_x(t) - p_x (t-s) ) ds \right)^2 dx.\\
\end{aligned}
\end{equation}
\begin{proof}
Taking the derivative of $\chi_1(t)$ defined in (\ref{3.1}), Using equations of (\ref{1.1}), and performing some integration by parts, we arrive at
\begin{equation*}
\begin{array}{ll}
\chi^\prime_1(t) &= \rho_1 || \phi_t(t)||^2_2+ \rho_2 || \psi_t(t)||^2_2
+  \rho_1 \int_0^L \phi \phi_{tt} ~dx  + \rho_2 \int_0^L \psi \psi_{tt}~dx\\
& = \rho_1 || \phi_t (t)||^2_2 + \rho_2 || \psi_t (t)||^2_2 - b || \psi_x(t)||^2_2- \kappa  h(t)|| p_x (t)||^2_2\\
&~~~ +~\kappa \int_0^L p_x(t) \int^\infty_0 g(s)(p_x(t)-p_x(t-s)) ds dx .\\
\end{array}
\end{equation*}
Applying Cauchy-Schwarz's, Young's inequalities and (\ref{2.15}) we obtain
\begin{equation*}
\chi^\prime_1(t) \leq  \rho_1 || \phi_t (t)||^2_2 + \rho_2 ||
\psi_t(t)||^2_2 - b || \psi_x(t) ||^2_2 - \kappa \biggl( h(t) -
\frac{\ell}{2} \biggr) || p_x(t)||^2_2 + \frac{\kappa}{ 2\ell} C_{\alpha}
(\mu \circ p_x)(t).
\end{equation*}
Note that $ h(t) - \frac{l}{2} \geq \frac{h(t)}{2}$ for all $t > 0 $, which proved (\ref{3.2}).\\\\
Deriving the functional $\chi_2 (t)$ set in (\ref{3.3}) and using the first equation of (\ref{1.1}) we get
\begin{equation} \label{3.5}
\begin{array}{ll}
\chi^\prime_2(t) &= - \rho_1 \int_0^L \phi_t \frac{\partial}{\partial t} \left[ \int_0^\infty g(s) ( p(t)-p(t-s) ) ds   \right]dx \\
&~~~- \rho_1 \int_0^L \phi_{tt} \int_0^\infty g(s) \left( p(t)-p(t-s) \right) ds dx  \\
& = - \rho_1 ( 1-h(t)) || \phi_t(t)||^2_2 \\
&~~~ - \rho_1 \int_0^L \phi_t \left[ \int_0^\infty g^\prime(s) \left(p(t)-p(t-s) \right)ds + (1-h)\widehat{\psi}_t (t) \right]dx \\
&~~~ +~ \kappa    \int_0^L \left( \int_0^\infty g(s) (  p_x(t) - p_x (t-s) ) ds \right)^2 dx \\
&~~~ + \kappa  h(t) \int_0^L p_x(t) \int_0^\infty  g(s) (p_x(t)-p_x(t-s)) ds dx.
\end{array}
\end{equation}
Applying again Cauchy-Schwarz's and Young's inequalities one has
\begin{multline} \label{3.6}
\left| \int_0^L \phi_t \int_0^\infty  g^{\prime}(s) (p(t)-p(t-s)) ds dx \right|  \\ \text{~~~~~~~~~~~~} \leq
\frac{g_0}{4} || \phi_t(t)||^2_2 + \frac{1}{g_0} \int_0^L \left( \int_0^\infty  g^{\prime}(s) (p(t)-p(t-s)) ds \right)^2 dx,
\end{multline}
\begin{equation} \label{3.7}
\left|  \int_0^L \phi_t \widehat{  \psi}_t(t) dx \right| \leq
\frac{1}{2} || \phi_t(t)||^2_2 + \frac{1}{2}
||  \widehat{\psi}_t(t)||^2_2 ,
\end{equation}
\begin{multline} \label{3.8}
 \left| \int_0^L p_x \int_0^\infty g(s)  (p_x(t)-p_x(t-s))dsdx \right| \\ \text{~~~~~~~~~~~~}\leq  (1- g_0)
\frac{\rho_1 L^2}{\rho_2} || p_x(t)||^2_2 +
\frac{\rho_2}{4(1- g_0)\rho_1 L^2}
 C_{\alpha} (\mu \circ p_x)(t).
\end{multline}
Replacing (\ref{3.6}-\ref{3.8}) in (\ref{3.5}) we obtain
\begin{equation}  \label{3.9}
\begin{array}{ll}
\chi^\prime_2(t) \leq & -\frac{\rho_1}{2} \Bigl( 1 - h(t)- \frac{g_0}{2} \Bigr)|| \phi_t(t)||^2_2 + \kappa  (1- g_0) \frac{\rho_1 L^2}{\rho_2} h(t) || p_x(t)||^2_2 \\
 & + \frac{\rho_1}{2} (1 -h(t))||  \widehat {\psi} _t(t)||^2_2 + \kappa \biggl( 1 + \frac{\rho_2}{4(1- g_0) \rho_1 L^2} \biggr) C_{\alpha} (\mu \circ p_x)(t) \\
 & + \frac{\rho_1}{g_0}  \int_0^L \left( \int_0^\infty g^{\prime}(s) (  p(t) - p (t-s) ) ds \right)^2 dx .
\end{array}
\end{equation}
Since, $ 1- h(t) = \int_0^t g(s) ds \geq \int^{t_0}_0 g(s) ds := g_0
> C_0,~ \forall  t \geq  t_0$, we recall that $1-h(t) \geq  g_0$,
and using standard computations, we conclude from (\ref{3.9}) that
the estimation (\ref{3.4}) satisfied.\\\\
Finally, by deriving $\chi_3(t)$ given in (\ref{3.10}) we obtain
\begin{equation} \label{3.12}
\begin{array}{ll}
\chi^\prime_3(t) = &  -\rho_2 || \psi_t(t)||^2_2
+ \kappa [h(t)]^2 || p_x(t)||^2_2 - \rho_2 \int_0^L \phi_{tt} \psi_x dx \\
& + 2 \kappa  h(t)   \int_0^L  p_x(t) \int_0^\infty g(s) (  p_x(t) - p_x (t-s) ) ds dx\\
& +\kappa  \int_0^L \left( \int_0^\infty g(s) (  p_x(t) - p_x (t-s) ) ds \right)^2 dx \\
& + \rho_2 \int_0^L \psi_t \left( g p_x - \int_0^\infty g^{\prime}(s) (  p_x(t) - p_x (t-s) ) ds \right) dx\\
&  + b \int_0^L \psi_{xx} \left( g p_x + \int_0^\infty g(s) (  p_x(t) - p_x (t-s) ) ds \right) dx.
\end{array}
\end{equation}
Integrating by parts the last term of (\ref{3.12}) and using first equation of (\ref{1.1}) we get
\begin{equation} \label{3.13}
\begin{array}{ll}
\chi^\prime_3(t) =& - \rho_2 || \psi_t(t)||^2_2
+ \kappa [h(t)]^2 || p_x(t) ||^2_2 + \biggl( \frac{b \rho_1}{\kappa} - \rho_2 \biggr) \int_0^L \phi_{tt} \psi_x dx\\
&+ 2 \kappa  h(t)   \int_0^L  p_x(t) \int_0^\infty g(s) (  p_x(t) - p_x (t-s) ) ds dx\\
&  +\kappa  \int_0^L \left( \int_0^\infty g(s) (  p_x(t) - p_x (t-s) ) ds \right)^2 dx \\
& + \rho_2 \int_0^L \psi_t \left( g p_x - \int_0^\infty g^{\prime}(s) (  p_x(t) - p_x (t-s) ) ds \right) dx.\\
\end{array}
\end{equation}
Using over again Cauchy-Schwarz's, Young's inequalities and (\ref{2.15}), and recalling that $g(t) \leq
g(0)$ for all $t > 0$, we deduce
\begin{multline} \label{3.14}
\left| \int_0^L \psi_t \left( g p_x - \int_0^\infty g^{\prime}(s) (  p_x(t) - p_x (t-s) ) ds \right) dx \right|  \\\text{~~~~~~} \leq \frac{1}{2} || \psi_t(t)||^2_2 + g(0) g(t)|| p_x||^2_2 + \int_0^L \left( \int_0^\infty g^{\prime}(s) (  p_x(t) - p_x (t-s) ) ds \right)^2 dx,
\end{multline}
\begin{equation} \label{3.15}
\left| h(t) \int_0^L p_x \int_0^\infty g(s) (  p_x(t) - p_x (t-s) ) ds dx
\right|  \leq \frac{[h(t)]^2}{2} || p_x(t)||^2_2 +\frac{1}{2}
C_{\alpha} (\mu \circ p_x)(t).
\end{equation}
Replacing (\ref{3.14}) and (\ref{3.15}) in (\ref{3.13}) and recalling the condition $\frac{\kappa}{\rho_1}=\frac{b}{\rho_2}$ we have
\begin{equation} \label{3.16}
\begin{array}{ll}
\chi^\prime_3(t) \leq & - \frac{\rho_2}{2}
|| \psi_t(t)||^2_2 + 2\kappa  [h(t)]^2 || p_x(t)||^2_2 + 2 \kappa  C_{\alpha} (\mu \circ p_x)(t)+\rho_2 g(0) g(t) || p_x(t)||^2_2 \\
&  + \rho_2 \int_0^L \left( \int_0^\infty g^{\prime}(s) (  p_x(t) - p_x (t-s) ) ds \right)^2 dx .
\end{array}
\end{equation}
Thus, regarding the equal wave speeds assumption $( \frac{\kappa}{\rho_1} = \frac{b}{\rho_2})$, we conclude from (\ref{3.16}) that (\ref{3.11}) satisfied.
\end{proof}
\end{lemma}

\begin{lemma} \label{lemma3.4}
Let $N$, $\eta_1$, $\eta_2 > 0$ be constants (to be determined later) and
\begin{equation} \label{3.17}
 L(t) := N E(t)+ \eta_1 I_1(t) + \eta_2 I_2(t) + I_3(t),
\end{equation}
satisfy, for some constant $C>0$, the estimate
\begin{equation}\label{3.18}
L^\prime (t) \leq - C E (t) + \frac{1}{4} (g \circ p_x)(t)
-4(1-\ell)|| p_x(t)||^2_2, ~~~ \forall  t \geq  t_0,
\end{equation}
\begin{proof}
Taking the derivative of $\chi (t)$ and using the estimates (\ref{3.2}), (\ref{3.4}), and (\ref{3.11}), we have, for all $t \geq  t_0$,
\begin{equation}  \label{3.19}
\begin{array}{ll}
\chi^\prime (t) \leq & - \rho_1 \Bigl( \eta_2 \frac{g_0}{4}-\eta_1 \Bigr) || \phi_t(t)||^2_2 - \eta_1 b|| \psi_x(t)||^2_2 - \biggl( \frac{\rho_2}{2} -\eta_1 \rho_2- \eta_2
\frac{\rho_1 L^2}{2} \biggr) || \psi_t(t)||^2_2\\
& - \kappa \biggl( \frac{\eta_1}{2} - \eta_2 (1-g_0) \frac{\rho_1 L^2}{\rho_2} - 2 h(t) \biggr) h(t) || p_x(t)||^2_2 -\left(\frac{N}{2}-C_1C_{\alpha}\right) (\mu \circ p_x)(t)\\
& -\left(\frac{Nk}{4}- \rho_2 g(0) g(t) \right)|| p_x(t)||^2_2 \\
&+ \Biggl( \eta_2 \frac{\rho_1 c^2_p}{g_0} + \rho_2 \Biggr)
\int_0^L \left( \int_0^\infty g^{\prime}(s) (  p_x(t) - p_x (t-s) ) ds \right)^2 dx,
\end{array}
\end{equation}where $C_1=  \kappa \biggl[ \frac{\eta_1}{2 l} + \eta_2 \biggl( 1 + \frac{\rho_2}{4(1-g_0) \rho_1 L^2} \biggr) + 2 \biggr]$.
\\ Now using  \textit{\textbf{(A1)}}, Remark \ref{new1} and Remark \ref{new2}, we notice that 
$$ g_0 > C_0 \geq \frac{64 \rho_1 L^2}{64 \rho_1 L^2+\rho_2}, $$ then $$ \frac{32 (1-g_0)}{g_0} < \frac{\rho_2}{2 \rho_1 L^2}.$$
So, it is possible to select $\eta_2$ such that
\begin{equation} \label{3.20}
\frac{32(1- g_0)}{g_0} < \eta_2 < \frac{\rho_2}{2\rho_1 L^2}.
\end{equation}
Similarly, using Remark's \ref{new1} and \ref{new2} one can show that 
$$ 8(1 -g_0) < \frac{1}{4} \min \{ \eta_2
g_0, 1 \}$$
which implies that we can select $\eta_1$ such that 
\begin{equation} \label{3.20a}
8(1 -g_0) < \eta_1 < \frac{1}{4} \min \{ \eta_2
g_0, 1 \} 
\end{equation}
From the choices in (\ref{3.20}) and (\ref{3.20a}) we observe that\\
$\bullet ~~~ \eta_2 \frac{g_0}{4}- \eta_1 > 0, \\
\bullet ~~~  \frac{\rho_2}{2}- \eta_1 \rho_2 -\eta_2
\frac{\rho_1 L^2}{2} > 0, \\
\bullet~~~  \frac{\eta_1}{2} - \eta_2 (1-g_0) \frac{\rho_1 L^2}{\rho_2} -2 h(t) > \frac{3}{2} (1- g_0) > 0~~~ \text{for all}~ t \geq  t_0.$\\
In this case, combining (\ref{3.20}) and (\ref{3.19}) we arrive at
\begin{equation} \label{3.21}
\begin{array}{ll}
\chi^\prime (t) \leq & - C E(t)  -\left(\frac{N}{2}-C_1C_{\alpha}\right) (\mu \circ p_x)(t)\\
& -\left(\frac{Nk}{4}- \rho_2 g(0) g(t) \right)|| p_x(t)||^2_2,
\end{array}
\end{equation}Now, using remark \ref{r1:4a}, there is $0<\alpha_0<1$ such that if $\alpha <\alpha_0$, then
\begin{equation}\label{con:r1:4a}
\alpha C_{\alpha}<\frac{1}{8C_1}.
\end{equation}Next, we choose $N$ large enough  so that
$$\frac{Nk}{4}- \rho_2 g(0) g(t)>4(1-\ell)\text{       and         }\alpha=\frac{1}{2N}<\alpha_0,$$
Therefore,  (\ref{3.21}) reduces to (\ref{3.18}).
\end{proof}
\end{lemma}

\begin{lemma}\label{lemma4a}
The functional
$$I_{4}(t)= \int_{0}^{t}h(t-s)\| p_x (t)\|_2^2dsdx,$$ satisfies, along the solution of (\ref{1.1}), the estimate

\begin{equation}\label{4a:1}
I_4^{\prime}(t)\le -\frac{1}{2}(g \circ  p_x)(t)+3(1-\ell)\| p_x
(t)\|_2^2dx+\frac{1}{2}\int_{t}^{+\infty} g(s)\vert \vert
 p_x \vert \vert_2^{2} ds.
\end{equation}
where $h(t)=\int_{t}^{+\infty}g(s)ds.$
\end{lemma}
\begin{proof}
As in \cite{Mahdi}, we have the following
\begin{equation}\label{n1}
\begin{aligned}
h'(t)=-g(t),~~~\int_0^{t}g(t-s) ds=\int_0^{t}g(s)
ds=\int_0^{\infty}g(s) ds-\int_t^{\infty}g(s) ds=h(0)-h(t).
\end{aligned}
\end{equation}
Now, direct differentiation of $I_4$, we have
\begin{equation}\label{n2}
\begin{aligned}
I_4'(t)&= h(0) \vert \vert p_x(t) \vert
\vert_2^{2}+\int_0^{t}h'(t-s) \vert \vert p_x(s) \vert
\vert_2^{2} ds\\
&=h(0) \vert \vert p_x(t) \vert \vert_2^{2}-\int_0^{t} g(t-s)
\vert
\vert p_x(s) \vert \vert_2^{2} ds\\
&=h(0) \vert \vert p_x(t) \vert \vert_2^{2}-\int_0^{t} g(t-s)
\vert
\vert p_x(s) -p_x(t) +p_x(t) \vert \vert_2^{2} ds\\
&=h(0) \vert \vert p_x(t)\vert \vert_2^{2} -\int_0^{t} g(t-s)
\vert \vert p_x(s) -p_x(t)\vert \vert_2^{2}-2 \int_0^{t}
g(t-s) \vert \vert p_x(s) -p_x(t)\vert \vert_2 \vert \vert
p_x(t) \vert \vert_2\\&-\int_0^{t} g(t-s) \vert
\vert p_x(t) \vert \vert_2^{2} ds\\
&=h(t) \vert \vert p_x(t)\vert \vert_2^{2}+2 \int_0^{t} g(t-s)
\vert \vert p_x(t) -p_x(s)\vert \vert_2 \vert \vert
p_x(t)\vert \vert_2-\int_0^{t}
g(t-s) \vert \vert p_x(t) \vert \vert_2^{2} ds\\
&\leq h(t) \vert \vert p_x(t)\vert \vert_2^{2}-\int_0^{t} g(t-s)
\vert \vert p_x(t) \vert \vert_2^{2} ds\\&+2(1-\ell)\vert \vert
p_x(t)\vert \vert_2^{2}+\frac{\int_0^{t}g(s)
ds}{2(1-\ell)}\int_0^{t}
g(t-s) \vert \vert p_x(t) -p_x(s)\vert \vert_2^2 ds\\
&\leq 3(1-\ell)\vert \vert p_x(t) \vert \vert_2^{2}-\int_0^{t}
g(t-s) \vert \vert p_x(t) -p_x(s)\vert \vert_2^2 ds+
\frac{1}{2}\int_0^{t} g(t-s) \vert \vert p_x(t) -p_x(s)\vert
\vert_2^2 ds
\end{aligned}
\end{equation}
\begin{equation}\label{n2}
\begin{aligned}
&\leq 3(1-\ell)\vert \vert p_x(t) \vert \vert_2^{2}-
\frac{1}{2}\int_0^{t} g(t-s) \vert \vert p_x(t) -p_x(s)\vert
\vert_2^2 ds\\
&\leq 3(1-\ell)\vert \vert p_x(t) \vert \vert_2^{2}-
\frac{1}{2}\int_0^{\infty} g(t-s) \vert \vert p_x(t)
-p_x(s)\vert \vert_2^2 ds+\frac{1}{2}\int_t^{\infty} g(t-s) \vert
\vert p_x(t) -p_x(s)\vert \vert_2^2 ds\\
&\leq 3(1-\ell)\vert \vert p_x(t) \vert \vert_2^{2}- \frac{1}{2}
\left(g \circ p_x \right)(t)+\frac{1}{2}\int_t^{\infty} g(t-s)
\vert \vert p_x(t) -p_x(s)\vert \vert_2^2 ds.
\end{aligned}
\end{equation}
Then \eqref{4a:1} is established.
\end{proof}
\begin{lemma}\label{r13:St}
Assume that
\begin{equation}\label{equalspeed}
   \frac{K}{\r_1}=\frac{b}{\r_2}.
\end{equation}
Then, the energy functional satisfies, for all $t\in \mathbb{R}^+$,
the following estimate
\begin{equation}\label{E1:r13:St}
\int_{0}^{t} E(s) ds < \tilde{m} f(t),
\end{equation}where $f(t)=1+\int_{0}^{t} h_0(s)ds$ and $h_0$ is defined in Lemma 3.2.
\end{lemma}
\begin{proof}
 let $F(t)=L(t)+I_{4}(t)$, then, we obtain, for all $t\in
\mathbb{R}_+$,
\begin{equation*}
\begin{aligned}
&L^{\prime}(t)\le -m E(t)+\frac{1}{2} \int_{t}^{+\infty} g(s) \vert
\vert p_x \vert\vert_2^{2} ds ,
\end{aligned}
\end{equation*}where $m$ is some positive constant. Therefore,
using \eqref{h}, we obtain
\begin{equation}\label{H}
\begin{aligned}
m\int_{0}^{t}E(s)ds&\le
F(0)-F(t)+\frac{M_0}{2}\int_{0}^{t}\int_{0}^{+\infty}
g(\tau+s)\left(1+\vert (p_x)_{0x}(s)\vert ds\right)^{2} d\tau
ds\\&\le F(0)+\frac{M_0}{2}\int_{0}^{t} h_0(s)ds.
\end{aligned}
\end{equation}
Hence, we  get
\begin{equation}\label{H}
\int_{0}^{t}E(s)ds\le \frac{F(0)}{m}+\frac{M_0}{2m}\int_{0}^{t}
h_0(s)ds \leq \tilde{m}\bigg(1+\int_{0}^{t} h_0(s)ds\bigg),
\end{equation}where $\tilde{m}=\max \big\{\frac{F(0)}{m},\frac{M_0}{2m}\big\}$.
\end{proof}

\begin{corollary}\label{corollary-1}
There exists $0< q_0 < 1$ such that, for all $t \ge 0$, we have the
following estimate:
\begin{equation}\label{E1c:L2:St0}
\int_{0}^{t}g(s)\vert\vert p_x(t)-p_x(t-s)\vert\vert_2^{2}ds\le
\frac{1}{q(t)} H^{-1}\left(\frac{q(t)\mu(t)}{\xi(t)}\right)
\end{equation}
where $H$ is defined earlier in Remark 2.1.
\begin{equation}\label{E3d:r13:St}
\mu(t):=-\int_{0}^{t}g'(s)\vert\vert p_x(t)-p_x(t-s)\vert
\vert_2^{2}ds\le -cE^{\prime}(t)
\end{equation}and
\begin{equation}\label{q}
    q(t):=\frac{q_0}{f(t)}.
\end{equation}
\end{corollary}
\begin{proof}
As in \cite{Guesmiarecent}, using \eqref{2.12} and \eqref{E1:r13:St}, we have
\begin{equation}\label{E:31Nc}
\begin{aligned}
& \int_0^t \vert\vert p_x(t)-p_x (t-s) \vert\vert_2^{2} ds
 \\& \le 2 \int_0^L \int_0^t \bigg( \vert p_x(t)\vert^{2}+ \vert p_x(t-s) \vert^{2} \bigg)ds dx\\
& \le \frac{4}{\kappa \gamma} \int_0^t \left(E(t) + E(t-s) \right) ds dx\\
&\le \frac{8}{\kappa \gamma} \int_0^{t} E(s)ds dx \leq
\frac{8}{\kappa \gamma} \tilde{m}f(t), ~~~\forall~~ t \in
\mathbb{R}_+.
\end{aligned}
\end{equation}Thanks to \eqref{E1:r13:St}, we have for all $t \ge
0$ and for  $0< q_0 < \min\big\{1, \frac{\kappa
\gamma}{8\tilde{m}}\big\}$, $q(t)<1$ and
$$q(t)\int_{0}^{t}\vert\vert p_x(t)-p_x(t-s)\vert \vert_2^{2}ds
<1.$$ So, the proof of \eqref{E1c:L2:St0} can be archived as the one
given in \cite{Mus-2017}.
\end{proof}
\section{A decay result for equal speeds of wave propagation}\label{sec4}
In this section, we state and prove a new general decay result in
the case of equal speeds of wave propagation \eqref{equalspeed}. As
in \cite{al2020general}, we  introduce the following  functions:
\begin{equation}\label{G1}
    G_{1}(t):=\int_{t}^{1}\frac{1}{s H^{\prime}(s)}ds,
\end{equation}
\begin{equation}\label{G234}
    G_{2}(t)=t H^{\prime}(t),\quad G_3(t)=t
    (H')^{-1}(t), \quad G_4(t)=G_3^{*}(t).
\end{equation}
Further,  we introduce the class $S$ of functions $\chi:
\mathbb{R}_+ \rightarrow \mathbb{R}_+^{*}$ satisfying for fixed
$c_1, c_2 > 0$ (should be selected carefully  in \eqref{4.36}):
\begin{equation}\label{xi} \chi\in
C^1(\mathbb{R}_+),~~~\chi \leq 1,~~~~\chi' \leq 0,
\end{equation}
and
\begin{equation}\label{dif}
   c_2 G_4\left[\frac{c}{d} q(t)
h_0(t)\right]\leq c_1  \left(
G_2\bigg(\frac{G_5(t)}{\chi(t)}\bigg)-\frac{G_2\left(G_5(t)\right)}{\chi(t)}\right),
\end{equation}
where $d >0$,  $c$ is a generic positive constant which may change
from line to line, $h_0$ and $q$ are defined in Lemma \ref{cor2}  and Corollary \ref{corollary-1} and
\begin{equation}\label{psi}
G_5(t)=G_1^{-1}\bigg(c_1\int_0^t
    \xi(s)ds\bigg),
\end{equation}
\begin{theorem}\label{main:th2017:1a}
Let $(\vp_0,\vp_1),(\psi_0,\psi_1)\in H^1_0(0,L)\times L^2(0,L)$.
Assume that the hypotheses (A.1) and (A.2) hold and
\begin{equation} \label{equal wave}
\frac{K}{\r_1}=\frac{b}{\r_2}.
\end{equation} then there exists a
strictly positive constant $C$ such that the solution of \ref{1.1}
satisfies, for all $t \geq 0$,
 \begin{equation}\label{decay1e}
E(t)\leq \frac{C G_5(t)}{\chi(t)q(t)}.
\end{equation}
\end{theorem}
\begin{remark}\cite{al2020general} According to the properties  of $H$ introduced
in $(A2)$, $G_2$ is convex increasing and defines a bijection from
$\mathbb{R}_+$ to $\mathbb{R}_+$, $G_1$ is decreasing defines a
bijection from $(0,1]$ to $\mathbb{R}_+$, and $G_3$ and $G_4$ are
convex and increasing functions on $(0,r]$. Then the  set $S$ is not
empty because it contains $\chi(s)= \varepsilon G_5(s)$ for any
$0 <  \varepsilon  \leq 1$ small enough. Indeed, \eqref{xi} is
satisfied (since \eqref{G1} and \eqref{psi}). On the other hand, we
have $q(t) h_0(t)$ is nonincreasing, $ 0 < G_5 \leq 1$, and $H'$ and
$G_4$ are increasing, then \eqref{dif} is satisfied if
$$ c_2 G_4 \bigg [ \frac{c}{d} q_0 h_0(0) \bigg] \leq  \frac{c_1}{\varepsilon} \bigg( H' \bigg(\frac{1 }{\varepsilon}
    \bigg)- H'(1)\bigg) $$ which holds, for $0 < \varepsilon \leq 1$ small enough, since $\lim_{t \rightarrow + \infty} H'(t) = + \infty$. But
with the choice $\chi= \varepsilon G_5$, \eqref{decay1e} (below)
does not lead to any stability estimate. The idea is to choose
$\chi$ satisfy \eqref{xi} and \eqref{dif} such that \eqref{decay1e}
gives the best possible decay rate for $E$.
\end{remark}
\begin{remark} \cite{Mahdi}
The stability estimate \eqref{decay1e} holds for any $\chi$
satisfying \eqref{xi} and \eqref{dif}. But \eqref{decay1e} does not
lead in general to the asymptotic stability $ \lim_{t\to \infty}
E(t)=0$ (like in case of the choice $\chi= \varepsilon G_5$
indicated in Remark 4.1, where \eqref{decay1e} becomes just an upper
bound estimate for $E$). The idea is to choose $\chi$ satisfying
\eqref{xi} and \eqref{dif} such that \eqref{decay1e} gives the best
possible decay rate for $E$. This choice can be done by taking
$\chi$ satisfying \eqref{xi} and \eqref{dif} such that the decay
rate of the function in the right hand side of \eqref{dif} has the
closet decay rate to the one of the function in the left hand side
of \eqref{dif}. So  such choice of $\chi$ can be seen from each
specific considered functions $g$ and $\psi_{0x}$ (see the
particular example considered below).
\end{remark}
\begin{proof} Now we start proof our main Theorem \ref{main:th2017:1a}\\
We combine  \eqref{h},  \eqref{3.18}, \eqref{E1c:L2:St0},  to have,
for some $m>0$ and for any $t\geq 0$, we have
\begin{equation}\label{e1s4b}
\begin{aligned}
L^{\prime}(t) \le -m E(t)+\frac{c}{q(t)}H^{-1}\left(\frac{q(t)
\mu(t)}{\xi(t)}\right)+ c h_0(t).
\end{aligned}
\end{equation}
Without loss of generality, one can assume that $E(0) > 0$. For
$\varepsilon_{0} < r $, let the functional $\mathcal{F}$ defined by
$$\mathcal{F}(t):=H^{\prime}\left(\varepsilon_{0}\frac{E(t)q(t)}{E(0)}\right)L(t),$$
which satisfies $\mathcal{F} \sim E$. By noting that
$H^{\prime\prime}\geq0,$ $q'\leq0$ and $E'\leq0$, we get
\begin{equation}\label{p6:F10b}
\begin{aligned}
&\mathcal{F}^{\prime}(t)=\varepsilon_{0}\frac{(qE)^{\prime}(t)}{E(0)}H^{\prime\prime}\left(\varepsilon_{0}\frac{E(t)q(t)}{E(0)}\right)L(t)+H^{\prime}\left(\varepsilon_{0}\frac{E(t)q(t)}{E(0)}\right)L^{\prime}(t)\\
&\hspace{0.3in}\le -m
E(t)H^{\prime}\left(\varepsilon_{0}\frac{E(t)q(t)}{E(0)}\right)+\frac{c}{q(t)}
H^{\prime}\left(\varepsilon_{0}\frac{E(t)q(t)}{E(0)}\right)H^{-1}\left(\frac{q(t)
\mu(t)}{\xi(t)}\right)\\&\hspace{0.3in}+c h_0(t)
H^{\prime}\left(\varepsilon_{0}\frac{E(t)q(t)}{E(0)}\right).
\end{aligned}
\end{equation}
Let $H^{*}$ be the convex conjugate of $H$ in the sense of Young
(see \cite{Yong}), then
\begin{equation}\label{p6:conj0b}
H^{*}(s)=s(H^{\prime})^{-1}(s)-H\left[(H^{\prime})^{-1}(s)\right],\hspace{0.1in}\text{if}\hspace{0.05in}s\in
(0,H^{\prime}(r)]
\end{equation}and satisfies the following generalized Young inequality
\begin{equation}\label{p6:young0b}
A B\le H^{*}(A)+H(B),\hspace{0.15in}\text{if}\hspace{0.05in}A\in
(0,H^{\prime}(r)],\hspace{0.05in}B\in(0,r].
\end{equation}So, with $A=H^{\prime}\left(\varepsilon_{0}\frac{E(t)q(t)}{E(0)}\right)$
 and $B=H^{-1}\left(\frac{q(t) \mu(t)}{\xi(t)}\right)$ and  using  \eqref{p6:F10b}-\eqref{p6:young0b}, we arrive at
\begin{equation}\label{E:m:xi10b}
\begin{aligned}
&\mathcal{F}^{\prime}(t)\le -m
E(t)H^{\prime}\left(\varepsilon_{0}\frac{E(t)q(t)}{E(0)}\right)+\frac{c}{q(t)}
H^{*}\left(H^{\prime}\left(\varepsilon_{0}\frac{E(t)q(t)}{E(0)}\right)\right)+c
\left(\frac{ \mu(t)}{\xi(t)}\right)\\&\hspace{0.3in}+c h_0(t) H^{\prime}\left(\varepsilon_{0}\frac{E(t)q(t)}{E(0)}\right)\\
&\hspace{0.3in}\le  -m
E(t)H^{\prime}\left(\varepsilon_{0}\frac{E(t)q(t)}{E(0)}\right)+c\varepsilon_{0}\frac{E(t)}{E(0)}H^{\prime}\left(\varepsilon_{0}\frac{E(t)q(t)}{E(0)}\right)+c
\left(\frac{ \mu(t)}{\xi(t)}\right)\\&\hspace{0.3in}+c h_0(t)
H^{\prime}\left(\varepsilon_{0}\frac{E(t)q(t)}{E(0)}\right).
\end{aligned}
\end{equation}
So, multiplying \eqref{E:m:xi10b} by $\xi(t)$ and using
\eqref{E3d:r13:St} and the fact that
$\varepsilon_{0}\frac{E(t)q(t)}{E(0)}<r$ gives
\begin{equation*}
\begin{aligned}
&\xi(t) \mathcal{F}^{\prime}(t)\le -m  \xi(t)
E(t)H^{\prime}\left(\varepsilon_{0}\frac{E(t)q(t)}{E(0)}\right)+c
\xi(t)\varepsilon_{0}\frac{E(t)}{E(0)}H^{\prime}\left(\varepsilon_{0}\frac{E(t)q(t)}{E(0)}\right)\\&\hspace{0.3in}+c
\mu(t)q(t)+c \xi(t)h_0(t)H^{\prime}\left(\varepsilon_{0}\frac{E(t)q(t)}{E(0)}\right)\\
&\hspace{0.3in}\le-\varepsilon_{0}(\frac{m  E(0)}{\varepsilon_{0}}
-c )
\xi(t)\frac{E(t)}{E(0)}H^{\prime}\left(\varepsilon_{0}\frac{E(t)q(t)}{E(0)}\right)-c
E^{\prime}(t)+c \xi(t)
h_0(t)H^{\prime}\left(\varepsilon_{0}\frac{E(t)q(t)}{E(0)}\right).
\end{aligned}
\end{equation*}Consequently, recalling the definition of $G_2$ and choosing $\varepsilon_{0}$ small enough so that $k=(\frac{m  E(0)}{\varepsilon_{0}}-c)>0$, we obtain, for all $t \in \mathbb{R}_+$,
\begin{equation}\label{p6:main40b}
\begin{aligned}
\mathcal{F}_{1}^{\prime}(t)&\le -k
\xi(t)\left(\frac{E(t)}{E(0)}\right)H^{\prime}\left(\varepsilon_{0}\frac{E(t)q(t)}{E(0)}\right)+c
\xi(t)
h_0(t)H^{\prime}\left(\varepsilon_{0}\frac{E(t)q(t)}{E(0)}\right)\\&=-k\frac{\xi(t)}{q(t)}
G_{2}\left(\frac{E(t)q(t)}{E(0)}\right)+c \xi(t)
h_0(t)H^{\prime}\left(\varepsilon_{0}\frac{E(t)q(t)}{E(0)}\right),
\end{aligned}
\end{equation}
where $\mathcal{F}_{1}=\xi \mathcal{F}+c E \sim E$ and satisfies for
some $\alpha_{1},\alpha_{2}>0.$
\begin{equation}\label{p6:equiv20b}
\alpha_{1}\mathcal{F}_1(t)\le E(t)\le \alpha_{2}\mathcal{F}_1(t).
\end{equation}
 Since
$G^{\prime}_{2}(t)=H^{\prime}(t)+t H^{\prime\prime}(t),$ then, using
the strict convexity of $H$ on $(0,r],$ we find that
$G_{2}^{\prime}(t), G_{2}(t)>0$ on $(0,r].$ Let $d>0$, use the
general Young inequality \eqref{p6:young0b} on the last term in
\eqref{p6:main40b} with
$A=H^{\prime}\left(\varepsilon_{0}\frac{E(t)q(t)}{E(0)}\right)$ and
$B=[\frac{c}{d}h_0 (t)]$,  to get
\begin{equation}\label{essa1a}
\begin{aligned}
c h_0
(t)H^{\prime}\left(\varepsilon_{0}\frac{E(t)q(t)}{E(0)}\right)&=\frac{d}{q(t)}
 \left[\frac{c}{d }q(t) h_0(t)\right] \bigg(H'
\left(\varepsilon_{0}\frac{E(t)q(t)}{E(0)}\right)\bigg)\\&\leq
\frac{d}{q(t)}
G_3\bigg(H'\left(\varepsilon_{0}\frac{E(t)q(t)}{E(0)}\right)\bigg)+\frac{d}{q(t)}G_3^*\left[\frac{c}{d
}q(t) h_0(t)\right]\\&\leq \frac{d}{q(t)}
\left(\varepsilon_{0}\frac{E(t)q(t)}{E(0)}\right)\left(H'\left(\varepsilon_{0}\frac{E(t)q(t)}{E(0)}\right)\right)+\frac{d}{q(t)}
G_4\left[\frac{c}{d }q(t) h_0(t)\right]\\&\leq \frac{d}{q(t)}
G_2\left(\varepsilon_{0}\frac{E(t)q(t)}{E(0)}\right)+\frac{d}{q(t)}
G_4\left[\frac{c}{d }q(t) h_0(t)\right].
\end{aligned}
\end{equation} Now, combining \eqref{p6:main40b} and
\eqref{essa1a} and choosing $d$ small enough so that $k_1=(k-d)>0$,
we arrive at
\begin{equation}\label{R'}
\begin{aligned}
   & \mathcal{F}_{1}^{\prime}(t) \leq - k \frac{\xi (t)}{q(t)} G_2 \left(\varepsilon_{0}\frac{E(t)q(t)}{E(0)}\right) + \frac{d \xi(t)}{q(t)}
G_2\left(\varepsilon_{0}\frac{E(t)q(t)}{E(0)}\right)+\frac{d
\xi(t)}{q(t)}
G_4\left[\frac{c}{d }q(t) h_0(t)\right]\\
&\hspace{0.3 in}\leq  - k_1 \frac{\xi (t)}{q(t)} G_2
\left(\varepsilon_{0}\frac{E(t)q(t)}{E(0)}\right)+\frac{d
\xi(t)}{q(t)}G_4\left[\frac{c}{d }q(t) h_0(t)\right].
\end{aligned}
\end{equation}Using the equivalent property in \eqref{p6:equiv20b} and since
 $G_2$ is increasing, we have, for some $d_{0}=\frac{\alpha_1}{E(0)}>0$,
\begin{equation*}
G_2 \left(\varepsilon_{0}\frac{E(t)q(t)}{E(0)}\right)\geq G_2
\bigg(d_{0} \mathcal{F}_{1}(t)q(t)\bigg).
\end{equation*}
Letting $\mathcal{F}_{2}(t):=d_{0} \mathcal{F}_{1}(t)q(t)$ and
recalling that $q'\leq0$, then we obtain for some constant $c_1=d_0
k_1 > 0$ and $c_2=d_0 d >0,$
\begin{equation}\label{4.36}
    \mathcal{F}_{2}'(t)\leq -c_1\xi (t)G_2(\mathcal{F}_{2}(t))+c_2
    \xi(t)G_4\left[\frac{c}{d }q(t) h_0(t)\right].
\end{equation}%
Since $d_{0}q(t)$  is non-increasing, then use of the equivalent
property $\mathcal{F}_{1}\sim E$ implies that there exists $b_0 > 0$
such that $\mathcal{F}_{2}(t)\geq b_{0} E(t)q(t)$. Let $t\in
\mathbb{R}_+$ and $\chi(t)$ satisfy \eqref{xi} and \eqref{dif}.\\If
\begin{equation}\label{decay000}
b_0q(t) E(t) \leq 2 \frac{G_5(t)}{\chi(t)},
\end{equation}
then, we have
\begin{equation}\label{decay1}
E(t) \leq \frac{2 }{b_0} \frac{G_5(t)}{\chi(t)q(t)}.
\end{equation}
If
\begin{equation}\label{decay000}
b_0q(t) E(t) > 2  \frac{G_5(t)}{\chi(t)},
\end{equation}
then, for any $0\leq s \leq t$, we get
\begin{equation}\label{decay000}
b_0q(s) E(s) > 2  \frac{G_5(t)}{\chi(t)},
\end{equation}
since, $q(t) E(t)$ is nonincreasing function. Therefore, we have for
any $0\leq s \leq t$,
\begin{equation}\label{psi1}
\mathcal{F}_{2}(s) > 2 \frac{G_5(t)}{\chi(t)}.
\end{equation}
Using \eqref{psi1}, the definition of $G_2$, the fact that $G_2$ is
convex  and $0< \chi\leq 1$, we have, for any $0 \leq s \leq t$ and
$0 < \epsilon_1 \leq 1$,
\begin{equation}\label{g2a}
\begin{aligned}
&G_2\bigg(\epsilon_1 \chi(s)\mathcal{F}_{2}(s)-\epsilon_1
G_5(s)\bigg) \leq \epsilon_1 \chi(s)G_2\bigg(\mathcal{F}_{2}(s)-
\frac{G_5(s)}{\chi(s)}\bigg)\\
&\leq \epsilon_1 \chi(s)\mathcal{F}_{2}(s)
H'\bigg(\mathcal{F}_{2}(s)- \frac{G_5(s)}{\chi(s)}\bigg)-\epsilon_1
\chi(s)\frac{G_5(s)}{\chi(s)}H'\bigg(\mathcal{F}_{2}(s)-
\frac{G_5(s)}{\chi(s)}\bigg)\\
&\leq \epsilon_1 \chi(s)\mathcal{F}_{2}(s)
H'\bigg(\mathcal{F}_{2}(s)\bigg)-\epsilon_1
\chi(s)\frac{G_5(s)}{\chi(s)}H'\bigg(\frac{G_5(s)}{\chi(s)}\bigg).
\end{aligned}
\end{equation}
Now, we let
\begin{equation}\label{F3}
    \mathcal{F}_{3}(t)=\epsilon_1 \chi(t)\mathcal{F}_{2}(t)-\epsilon_1
G_5(t),
\end{equation} where $\epsilon_1$ is small enough so that $\mathcal{F}_{3}(0)\leq 1$. Then \eqref{g2a}
becomes, for any $0 \leq s \leq t$,
\begin{equation}\label{4.37}
\begin{aligned}
&G_2\bigg(\mathcal{F}_{3}(s)\bigg) \leq \epsilon_1
\chi(s)G_2\bigg(\mathcal{F}_{2}(s)\bigg)- \epsilon_1
\chi(s)G_2\bigg(\frac{G_5(s)}{\chi(s)}\bigg).
\end{aligned}
\end{equation}
Further, we have
\begin{equation}\label{g2}
\begin{aligned}
\mathcal{F}'_{3}(t)=\epsilon_1 \chi'(t)\mathcal{F}_{2}(t)
+\epsilon_1\chi(t)\mathcal{F}'_{2}(t)-\epsilon_1 G_5'(t).
\end{aligned}
\end{equation}
Since $\chi' \leq 0$ and using \eqref{4.36}, then for any $0 \leq s
\leq t$, $0 < \epsilon_1 \leq 1$,  we obtain
\begin{equation}\label{g2}
\begin{aligned}
&\mathcal{F}'_{3}(t)\leq
\epsilon_1\chi(t)\mathcal{F}'_{2}(t)-\epsilon_1 G_5'(t)\\
&\leq -c_1\epsilon_1 \xi (t)
\chi(t)G_2(\mathcal{F}_{2}(t))+c_2\epsilon_1 \xi (t)
\chi(t)G_4\left[\frac{c}{d }q(t) h_0(t)\right]-\epsilon_1G_5'(t).
\end{aligned}
\end{equation}
Then, using \eqref{dif} and \eqref{4.37}, we get
\begin{equation}\label{f3}
\begin{aligned}
&\mathcal{F}'_{3}(t) \leq -c_1 \xi (t)
G_2(\mathcal{F}_{3}(t))+c_2\epsilon_1 \xi
(t)\chi(t)G_4\left[\frac{c}{d }q(t) h_0(t)\right]\\&-c_1\epsilon_1
\xi
(t)\chi(t)G_2\bigg(\frac{G_5(t)}{\chi(t)}\bigg)-\epsilon_1G_5'(t).
\end{aligned}
\end{equation}
From the definition of $G_1$ and $G_5$, we have
\begin{equation*}
    G_1\left(G_5(s)\right)=c_1\int_0^s \xi(\tau) d\tau,
\end{equation*}
hence,
\begin{equation}\label{psi'}
 G_5'(s)=-c_1 \xi(s) G_2\left(G_5(s)\right).
\end{equation}
Now, we have
\begin{equation}\label{coll}
\begin{aligned}
&c_2\epsilon_1 \xi (t)\chi(t)G_4\left[\frac{c}{d }q(t)
h_0(t)\right]-c_1\epsilon_1 \xi
(t)\chi(t)G_2\bigg(\frac{G_5(t)}{\chi(t)}\bigg)-\epsilon_1G_5'(t)\\&=c_2\epsilon_1
\xi (t)\chi(t)G_4\left[\frac{c}{d }q(t) h_0(t)\right]-c_1\epsilon_1
\xi
(t)\chi(t)G_2\bigg(\frac{G_5(t)}{\chi(t)}\bigg)+c_1\epsilon_1\xi(t)G_2\left(G_5(t)\right)\\&=\epsilon_1
\xi (t)\chi(t)\bigg(c_2 G_4\left[\frac{c}{d }q(t)
h_0(t)\right]-c_1G_2\bigg(\frac{G_5(t)}{\chi(t)}\bigg)+c_1\frac{G_2\left(G_5(t)\right)}{\chi(t)}\bigg).
\end{aligned}
\end{equation}
Then, according to \eqref{dif}, we get
$$\epsilon_1 \xi
(t)\chi(t)\bigg(c_2G_4\left[\frac{c}{d }q(t)
h_0(t)\right]-c_1G_2\bigg(\frac{G_5(t)}{\chi(t)}\bigg)-c_1\frac{G_2\left(G_5(t)\right)}{\chi(t)}\bigg)\leq
0 $$ Then \eqref{f3} gives
\begin{equation}\label{f30}
\begin{aligned}
&\mathcal{F}'_{3}(t) \leq -c_1 \xi (t) G_2(\mathcal{F}_{3}(t)).
\end{aligned}
\end{equation}
Thus  from \eqref{f30} and the definition of $G_1$ and $G_2$ in
\eqref{G1} and \eqref{G234}, we obtain
\begin{equation}\label{end1a}
    \bigg(G_1\left(\mathcal{F}_{3}(t)\right)\bigg)'\geq c_1
    \xi(t).
\end{equation}
Integrating  \eqref{end1a} over $[0, t]$, we get
\begin{equation}\label{end2}
 G_1\left(\mathcal{F}_{3}(t)\right)\geq c_1\int_0^t
    \xi(s)ds+G_1\left(\mathcal{F}_{3}(0)\right).
\end{equation}
Since $G_1$ is decreasing, $\mathcal{F}_{3}(0)\leq 1$ and
$G_1(1)=0$, then

\begin{equation}\label{end3}
\mathcal{F}_{3}(t)\leq G_1^{-1}\bigg(c_1\int_0^t
    \xi(s)ds\bigg)=G_5(t).
\end{equation}
Recalling that $ \mathcal{F}_{3}(t)=\epsilon_1
\chi(t)\mathcal{F}_{2}(t)-\epsilon_1 G_5(t)$, we have

\begin{equation}\label{end4}
\mathcal{F}_{2}(t)\leq \frac{(1+\epsilon_1)}{\epsilon_1}
\frac{G_5(t)}{\chi(t)},
\end{equation}
Similarly, recall that $\mathcal{F}_{2}(t):=d_{0}
\mathcal{F}_{1}(t)q(t)$, then
\begin{equation}\label{end5}
\mathcal{F}_{1}(t)\leq \frac{(1+\epsilon_1)}{d_0\epsilon_1}
\frac{G_5(t)}{\chi(t)q(t)},
\end{equation}
Since $\mathcal{F}_{1}\sim E$, then for some $b>0$, we have
$E(t)\leq b \mathcal{F}_{1}$; which gives
\begin{equation}\label{decay2}
E(t)\leq \frac{b(1+\epsilon_1)}{d_0\epsilon_1}
\frac{G_5(t)}{\chi(t)q(t)},
\end{equation}
From \eqref{decay1} and \eqref{decay2}, we obtain the following
estimate
\begin{equation}\label{decay3}
E(t)\leq c_3 \bigg(\frac{G_5(t)}{\chi(t)q(t)}\bigg),
\end{equation}
where $c_3=\max\{\frac{2}{b_0},
\frac{b(1+\epsilon_1)}{d_0\epsilon_1}\}.$ This complete the proof of
Theorem 4.1.
\end{proof}

\textbf{Example 4.1 \cite{al2020general}, \cite{Guesmiarecent}:} Let
$g(t)=\frac{a}{(1+t)^\nu}$, where $\nu
>1$ and $0<a<\nu-1$ so that $(A1)$ is satisfied. In this case
$\xi(t)=\nu a^{\frac{-1}{\nu}}$ and $G(t)=t^{\frac{\nu+1}{\nu}}$.
Then, there exist positive constants $a_i (i = 0,...,3)$ depending
only on $a, \nu$ such that
\begin{equation}\label{gs}\begin{aligned}
&G_4(t)=a_0 t^{\frac{\nu+1}{\nu}},~~~G_2(t)=a_1
t^{\frac{\nu+1}{\nu}},~~~~G_1(t)=a_2
(t^{\frac{-1}{\nu}}-1),\\
& G_5(t)=(a_4 t+1)^{-\nu},~~G_3(t)=a_3 t^{\nu +1}.
\end{aligned}
\end{equation}
We  will discuss two cases:\\Case 1: if
\begin{equation}\label{r}
    m_0 (1+t)^{r} \leq 1+\vert \vert {{p_x}}_{0x} \vert \vert^2 \leq m_1
(1+t)^{r}
\end{equation}
where $0 < r < \nu-1$ and $m_0 , m_1 > 0$, then we have, for some
positive constants $a_i (i = 5,...,8)$ depending only on $a, \nu,
m_0, m_1, r$,  the following:
\begin{equation}\label{gs0}\begin{aligned}
& a_5 (1+t)^{-\nu+1+r}\leq h_0(t) \leq a_6 (1+t)^{-\nu+1+r},
\end{aligned}
\end{equation}
\begin{equation}\label{gs10}\begin{aligned}
&\frac{q_0}{q(t)} \geq a_7\left\{%
\begin{array}{ll}
    1+\ln(1+t), & \hbox{$\nu-r=2$;} \\
    2, & \hbox{$\nu -r> 2$;} \\
    (1+t)^{-\nu+r+2}, & \hbox{$1< \nu -r <2$ .} \\
\end{array}%
\right.
\end{aligned}\end{equation}
\begin{equation}\label{gs20}\begin{aligned}
&\frac{q_0}{q(t)} \leq a_8\left\{%
\begin{array}{ll}
    1+\ln(1+t), & \hbox{$\nu-r=2$;} \\
    2, & \hbox{$\nu -r> 2$;} \\
    (1+t)^{-\nu+r+2}, & \hbox{$1< \nu -r <2$ .} \\
\end{array}%
\right.
\end{aligned}
\end{equation}We notice that condition
\eqref{dif} is satisfied if
\begin{equation}\label{new}
    (t+1)^{\nu} q(t) h_0(t) \chi(t) \leq a_9 \bigg(1- (\chi)^{\frac{1}{\nu}}   \bigg)^{\frac{\nu}{\nu
    +1}}.
\end{equation}
where $a_9 > 0$  depending on $a, \nu, c_1$ and $c_2$. Choosing
$\chi(t)$ as the following
\begin{equation}\label{gs3}\begin{aligned}
&\chi(t) = \lambda\left\{%
\begin{array}{ll}
     (1+t)^{-p},~~~~~ p=r+1 & \hbox{$\nu-r \geq 2$;} \\
   (1+t)^{-p},~~~~~ p=\nu-1, & \hbox{$1< \nu -r <2$ .} \\
\end{array}%
\right.
\end{aligned}
\end{equation}
so that \eqref{xi} is valid. Moreover, using \eqref{gs0} and
\eqref{gs10}, we see that \eqref{new} is satisfied if $0 < \lambda
\leq 1$ is small enough, and then \eqref{dif} is satisfied. Hence
\eqref{decay1e} and \eqref{gs20} imply that, for any $t \in
\mathbb{R}_+$
\begin{equation}\label{gsw1}\begin{aligned}
&E(t) \leq a_{10}\left\{%
\begin{array}{ll}
    \bigg(1+\ln(1+t)\bigg)(1+t)^{-(\nu-r-1)}, & \hbox{$\nu-r=2$;} \\
  (1+t)^{-(\nu-r-1)}, & \hbox{$\nu-r>2$;} \\
    (1+t)^{-(\nu-r-1)}, & \hbox{ $1< \nu-r <2$ .} \\
\end{array}%
\right.
\end{aligned}
\end{equation}
Thus,  the  estimate \eqref{gsw1} gives $\lim_{t\rightarrow +\infty}
E(t)=0$.\\
Case 2: if $m_0 \leq 1+\vert \vert {{p_x}}_{0x} \vert \vert^2 \leq
m_1$. That is $r = 0$ in (\ref{gsw1}).\\\\


\section*{Acknowledgements}
The author thanks King Fahd University of Petroleum and Minerals (KFUPM) for the continuous support. The author also thanks the referee for her/his very careful reading and valuable comments. This work is funded by KFUPM under project \# SR181024.




\end{document}